\documentclass[11pt]{article}
\usepackage{latexsym,amssymb,amsmath,amscd,amsfonts}
\usepackage{graphicx}
\usepackage{mathrsfs}
\usepackage{xcolor}
\usepackage[dvipsnames]{pstricks}
\usepackage{microtype}
\usepackage{needspace}
\usepackage[hidelinks]{hyperref}
\allowdisplaybreaks
\newtheorem{thm}{Theorem}[section]

\newtheorem{pro}[thm]{Proposition}

\newtheorem{defi}[thm]{Definition}

\newtheorem{core}[thm]{Corollary}
\newtheorem{exam}[thm]{\bf Example}
\def\pf{\noindent{\it Proof.} \nopagebreak[4]}
\def\qed{\nopagebreak\hfill{\rule{4pt}{7pt}}
	\medbreak}

\numberwithin{equation}{section}

\makeatletter
\def\ExtendSymbol#1#2#3#4#5{\ext@arrow 0099{\arrowfill@#1#2#3}{#4}{#5}}
\makeatother

\title{$k$-Pairing: A Generalization of the Partition Pairing Theorems}
\author{
	Xiaorui Niu\raisebox{5pt}{\scriptsize 1}
	and
	Diane Y. H. Shi\raisebox{5pt}{\scriptsize 2}
}
\date{
	\vspace{15pt}
	\raisebox{5pt}{\scriptsize 1\,}xiaorui\_niu@tju.edu.cn\\
	School of Mathematics\\
	Tianjin University, Tianjin 300072, P.~R.~China\\[10pt]
	\raisebox{5pt}{\scriptsize 2\,}shiyahui@tju.edu.cn\\
	School of Mathematics\\
	Tianjin University, Tianjin 300072, P.~R.~China
}

\begin{document}

\maketitle

\noindent {\bf Abstract.}
We extend the partition pairing theory of Andrews and Dastidar by replacing pairs with groups of \(k\) equal parts. Two weight-preserving bijections give combinatorial interpretations of the joint pairing index--width distribution and the negative-rank enumerations. The first combines conjugation with the Stockhofe--Keith correspondence and sends the \(k\)-pairing index and width to the number of parts and the largest part, respectively. Consequently, their joint distribution is independent of \(k\) and is given by a Gaussian polynomial. We also obtain finite refinements that record the residual statistics.
The second bijection sends simply \(k\)-paired partitions of negative \(k\)-pairing rank to \(k\)-regular partitions with marked internal gaps. In a fixed nonzero residue class, sign cancellation leaves rectangular partitions, while gap markings correspond to overlining choices with the smallest part not overlined. This gives direct combinatorial explanations of the divisor counts and the factor \(1/2\) in the overpartition enumeration, extending the odd-divisor and odd-overpartition results of Andrews and Dastidar.
Finally, motivated by the diagonal pairing of Andrews and Dastidar, we extend the construction to ordered tuples of \(k\) Young diagrams. Applying the two-wing transfer operation to pairs of component diagrams defines an equivalence relation on these tuples. We show that two tuples are equivalent if and only if they have the same cell-multiplicity function, and that every equivalence class contains a unique representative whose component diagrams are nested. We also determine the cardinalities of the equivalence classes and identify the nested representatives with plane partitions of rectangular shape.



\section{Introduction}
\label{sec:introduction}

We begin by recalling the basic notation for partitions and the partition pairing theory of
Andrews and Dastidar. Throughout, a partition is a finite weakly decreasing sequence of positive integers
\[
\lambda=(\lambda_1,\ldots,\lambda_\ell).
\]
We write \(|\lambda|=\lambda_1+\cdots+\lambda_\ell\) for its size, \(\ell(\lambda)=\ell\) for its
length, \(m_j(\lambda)\) for the multiplicity of the part \(j\), and \(\lambda'\) for its conjugate.
We also write \(\varnothing\) for the empty partition and use \(\sqcup\) for multiset union.

Andrews and Dastidar~\cite{and26} introduced partition pairing by pairing equal parts as far as
possible. If \(r(\lambda)\) is the largest paired part, with \(r(\lambda)=0\) when there is no pair,
and the unpaired parts are \(u_1>\cdots>u_s\), they define
\[
s(\lambda)=u_1-u_2+u_3-\cdots.
\]
The pairing index and pairing rank are
\[
T(\lambda)=r(\lambda)+s(\lambda),
\qquad
\rho(\lambda)=r(\lambda)-s(\lambda),
\]
respectively. They proved that \((r,s)\) is jointly equidistributed with the numbers of even and
odd parts. Their pairing width \(W\) further yields the joint equidistribution of \((T,W)\) with
the number of parts and the largest part, with the corresponding generating function given by
a Gaussian polynomial.

We replace pairs by groups of \(k\) equal parts. Writing
\(m_j(\lambda)=ka_j+\varepsilon_j\), where \(0\leq\varepsilon_j<k\), gives a unique decomposition
\(\lambda\longleftrightarrow(\pi,\delta)\): the multiplicities of \(\pi\) are \(a_j\), and those
of \(\delta\) are \(\varepsilon_j\). For \(k=2\), the partition \(\delta\) has distinct
parts. For general \(k\), we use conjugation and the Stockhofe--Keith bijection, together with
the residue and largest-part refinements of Xiong and Keith~\cite{xio19} and Fu, Tang, and
Yee~\cite{fu19}.

Let \(\mathcal D_k\) be the set of partitions in which every multiplicity is less than \(k\), and
let \(\mathcal O_k\) be the set of \(k\)-regular partitions, whose parts are not divisible by
\(k\). Denote the resulting weight-preserving bijection by
\(\Theta_k:\mathcal D_k\longrightarrow\mathcal O_k\). We define
\[
  \Phi_k(\lambda)=k\pi'\sqcup\Theta_k(\delta).
\]
This gives a weight-preserving bijection on ordinary partitions and yields
\[
  \bigl(T_k(\lambda),W_k(\lambda)\bigr)
  =\bigl(\ell(\Phi_k(\lambda)),\Phi_k(\lambda)_1\bigr).
\]
More precisely, the numbers of parts of \(\Phi_k(\lambda)\) in the residue classes modulo \(k\)
record the separate constituents of \(T_k\), while its largest part is \(W_k\). Consequently,
if \(N_{t,b}^{(k)}(n)\) counts partitions of \(n\) with \(T_k=t\) and \(W_k=b\), then
\[
  \sum_{n\geq0}N_{t,b}^{(k)}(n)q^n
  =q^{t+b-1}\begin{bmatrix}t+b-2\\t-1\end{bmatrix}_q
  \qquad(t,b\geq1).
\]
Thus the joint distribution of \(T_k\) and \(W_k\) is independent of \(k\), and the
corresponding Gaussian-polynomial identity in~\cite{and26} is recovered when \(k=2\).

We next consider simply \(k\)-paired partitions, in which no part occurs more than \(2k-1\)
times. For negative \(k\)-pairing rank, we construct a weight-preserving bijection
\[
  \Psi_k(\lambda)=\bigl(\Theta_k(\delta)+kA',A\bigr),
\]
where \(A\) is the set of part sizes occurring at least \(k\) times, regarded as a distinct
partition when taking its conjugate. Each group of \(k\) equal parts gives a mark on an internal
gap of the output partition. The bijection \(\Psi_k\) gives a combinatorial interpretation of the
signed and unsigned enumerations for negative \(k\)-pairing rank and extends the results of
Andrews and Dastidar~\cite{and26} on odd divisors and overpartitions into odd parts.
Fix \(1\leq a<k\) and require every part of \(\Theta_k(\delta)\) to be congruent to \(a\pmod k\).
Let \(E_{k,a}^{\mathrm{sp}}(n)\) and \(O_{k,a}^{\mathrm{sp}}(n)\) count simply \(k\)-paired
partitions of \(n\) satisfying this condition and \(\rho_k(\lambda)<0\), with an even and an odd
number of groups of \(k\) equal parts, respectively. Put
\(M_{k,a}(n)=E_{k,a}^{\mathrm{sp}}(n)+O_{k,a}^{\mathrm{sp}}(n)\).
For \(n\geq1\), we obtain
\[
\begin{aligned}
 E_{k,a}^{\mathrm{sp}}(n)-O_{k,a}^{\mathrm{sp}}(n)
   &=\#\{d\mid n:d\equiv a\pmod k\},\\
 M_{k,a}(n)&=\tfrac12\,\overline p_{a,k}(n),
\end{aligned}
\]
where \(\overline p_{a,k}(n)\) counts overpartitions with all parts congruent to \(a\pmod k\);
see~\cite{cor04}. We also obtain generating functions for fixed values of \(\rho_k\).
Under the same restriction on the parts of \(\Theta_k(\delta)\), we derive a recurrence relating
\(A_d^{(k,a)}(q)\) and \(A_{d+1}^{(k,a)}(q)\).

Finally, we extend the diagonal pairing of~\cite{and26} to ordered tuples
\((d;\alpha^{(0)},\ldots,\alpha^{(k-1)})\), where each component partition has at most \(d\)
parts. The size of a tuple is \(d^2+\sum_s|\alpha^{(s)}|\).
Choose two component diagrams and transfer an edge-connected component of their symmetric
difference from one diagram to the other, leaving the remaining components unchanged.
Two tuples are equivalent if one can be obtained from the other by a sequence of these
operations. For fixed \(d\), we prove that they are equivalent if and only if each cell occurs
in the same number of component diagrams. Put \(h(x)=\#\{s:x\in\alpha^{(s)}\}\) and
\(U_r=\{x:h(x)\geq r\}\). Every equivalence class has a unique representative whose component
diagrams are weakly decreasing under inclusion, namely \((d;U_1,\ldots,U_k)\).

For \(k=2\), each edge-connected component of \(U_1\setminus U_2\) can be assigned independently
to either diagram, giving a class of size \(2^\kappa\), where \(\kappa\) is the number of these
components. For \(k\geq3\), the assigned subsets must satisfy the inclusions in
\eqref{eq:sec4-15}. The number \(\Omega_k(U_1,\ldots,U_k)\) defined there gives the class
cardinality, for which we also prove a product formula. For fixed \(d\), the representatives
\((d;U_1,\ldots,U_k)\) correspond to plane partitions whose shape is contained in a
\(k\times d\) rectangle. Cyclic shifts of the component partitions give root-of-unity
evaluations that count fixed tuples. Finally, for \(m\geq2\), the map
\(\lambda\mapsto\alpha=(\Phi_m(\lambda))'\) satisfies
\(\ell(\alpha)=W_m(\lambda)\) and \(\alpha_1=T_m(\lambda)\).

This paper is  organized as follows. In Section~\ref{sec:k-pairing} we shall state the bijection \(\Phi_k\) and the resulting generating
functions, including formulas with bounds on \(W_k(\lambda)\) or \(\lambda_1\).
Section~\ref{sec:negative-rank} treats negative rank and the bijection \(\Psi_k\).
Section~\ref{sec:partition-tuples} studies equivalence classes of tuples of partitions and their
connection with pairing width.

\section{\texorpdfstring{$k$}{k}-Pairing and Joint Distributions}
\label{sec:k-pairing}

Let \(\mathcal P\) denote the set of all partitions, including the empty partition. The largest part
of the empty partition is taken to be zero, and partitions are padded with zeros when needed.
We employ the standard \(q\)-series notation. For \(N\geq0\), we write
\[
(a;q)_N=\prod_{j=0}^{N-1}(1-aq^j),
\qquad
(a;q)_\infty=\prod_{j\geq0}(1-aq^j),
\]
and use the Gaussian coefficient
\[
\begin{bmatrix}
A\\ B
\end{bmatrix}_q
=
\begin{cases}
\dfrac{(q;q)_A}{(q;q)_B(q;q)_{A-B}},
&0\leq B\leq A,\\[3mm]
0,&\text{otherwise}.
\end{cases}
\]

We also use the usual \(q\)-multinomial coefficient
\[
\begin{bmatrix}
N\\n_0,n_1,\ldots,n_s
\end{bmatrix}_q
=
\frac{(q;q)_N}
{(q;q)_{n_0}(q;q)_{n_1}\cdots(q;q)_{n_s}},
\]
whenever \(n_0+\cdots+n_s=N\).

\subsection{The Stockhofe--Keith bijection and alternating sums}
\label{subsec:stockhofe-keith}

Fix an integer \(k\geq2\). Let \(\mathcal D_k\) be the set of partitions in which every multiplicity
is less than \(k\), and let \(\mathcal O_k\) be the set of partitions with no part divisible by
\(k\).

Following the \(m\)-alternating sum statistics used by Xiong and Keith~\cite{xio19} and by Fu, Tang, and Yee~\cite{fu19}, we introduce their \(m=k\) specialization in the present notation. For \(\delta=(\delta_1,\delta_2,\ldots)\in\mathcal D_k\) and \(1\leq i\leq k\), define
\begin{equation}
\label{eq:sec2-1}
\Sigma_i^{(k)}(\delta)
=
\sum_{t\geq0}
\bigl(
\delta_{tk+i}-\delta_{tk+i+1}
\bigr).
\end{equation}

Their combinatorial meaning becomes transparent after conjugation. Put \(\tau=\delta'\). For each \(j\geq1\), \(\delta_j-\delta_{j+1}\) is the number of columns of the Ferrers diagram of \(\delta\) having height exactly \(j\). Equivalently, \(\delta_j-\delta_{j+1} = m_j(\tau)\). Thus the summands occurring in \eqref{eq:sec2-1} count columns according to their heights modulo \(k\). We therefore have
\[
\Sigma_i^{(k)}(\delta)
=
\#\{
\text{columns of }\delta:
\text{height}\equiv i\pmod k
\}
\]
or, equivalently,
\[
\Sigma_i^{(k)}(\delta)
=
\#\{
\text{positive parts of }\tau:
\tau_j\equiv i\pmod k
\}.
\]

Here \(i=k\) represents the residue class \(0\pmod k\).

We put \(s_k(\delta)=\sum_{i=1}^{k-1}\Sigma_i^{(k)}(\delta)\) and
\(C_k(\delta)=\Sigma_k^{(k)}(\delta)=\delta_1-s_k(\delta)\). The last equality follows from
\[
\sum_{i=1}^{k}\Sigma_i^{(k)}(\delta)
=\sum_{j\geq1}(\delta_j-\delta_{j+1})
=(\delta_1-\delta_2)+(\delta_2-\delta_3)+\cdots=\delta_1.
\]

We next describe the map \(\Theta_k\) used below.

A partition \(\tau=(\tau_1,\tau_2,\ldots,\tau_\ell)\) is called \(k\)-flat if
\(0\leq\tau_j-\tau_{j+1}<k \quad(1\leq j<\ell)\) and \(0<\tau_\ell<k\). The empty partition is also
regarded as \(k\)-flat.

Since \(\tau_j-\tau_{j+1} = m_j(\delta)\), the condition \(\delta\in\mathcal D_k\) gives
\(0\leq\tau_j-\tau_{j+1}\leq k-1\). If we regard the part following the smallest positive part as
\(0\), the same inequality also yields \(\tau_\ell<k\). Hence conjugation gives a weight-preserving
bijection between \(\mathcal D_k\) and the set of \(k\)-flat partitions.

We now recall the Stockhofe--Keith bijection in the form needed here;
see \cite[Section~3]{xio19} and \cite[Section~2.1]{fu19}.
For \(\delta\in\mathcal D_k\), put \(\tau=\delta'\), and let
\(\gamma=\tau\) denote the current partition.
Start with an empty list. Whenever weight \(ks\) is removed,
record the integer \(s\).

\smallskip
\noindent\textbf{Step 1.}
Delete any occurrence of a part \(ak\) whose removal leaves
\(\gamma\) \(k\)-flat, and record \(a\).
Continue until no further deletion of this type is possible.
The order of these deletions does not matter.

\smallskip
\noindent\textbf{Step 2.}
If a multiple of \(k\) remains, choose the largest one.
Suppose that it is the \(i\)-th part, of size \(ak\).
Replace the current partition by
\begin{equation}
\label{eq:sec2-8}
\begin{aligned}
&(\gamma_1,\ldots,\gamma_{i-1},ak,\gamma_{i+1},\ldots)\\
&\qquad\longmapsto
(\gamma_1-k,\ldots,\gamma_{i-1}-k,\gamma_{i+1},\ldots),
\end{aligned}
\end{equation}
and record \(a+i-1\), since the removed weight is
\(ak+k(i-1)\).
Renumber the parts and repeat until no multiple of \(k\) remains.

To see why the subtraction is needed, note that after Step~1
a remaining multiple \(ak\) has neighbouring parts \(ak+r\) and
\((a-1)k+s\), where \(1\le s\le r\le k-1\).
Direct deletion would create a gap \(k+r-s\ge k\).
In \eqref{eq:sec2-8}, the two newly adjacent parts instead differ
by \(r-s<k\). All other adjacent differences and the smallest
part are unchanged, so \(k\)-flatness is preserved.

\smallskip
\noindent\textbf{Step 3.}
Let \(\beta\) be the remaining partition, which is both
\(k\)-flat and \(k\)-regular.
Arrange the recorded integers in weakly decreasing order
to form \(\sigma\), and define
\begin{equation}
\label{eq:sec2-9}
\Theta_k(\delta)=\nu=\beta+k\sigma',
\end{equation}
where the addition is componentwise, with zero padding.

We check that this addition creates no new positive parts.
Let \(L\) be the number of parts of \(\tau\) not divisible by \(k\).
These parts are never deleted, so \(\ell(\beta)=L\).
Whenever a current part equals \(ak\), flatness forces a part
in each interval \((bk,(b+1)k)\), where \(0\le b<a\);
otherwise an adjacent difference, including the difference
from the smallest part to zero, would be at least \(k\).
Thus \(a\le L\) in Step~1.
In Step~2, the \(i-1\) preceding parts are all nonzero modulo \(k\),
and the \(a\) intervals lie below \(ak\), giving \(a+i-1\le L\).
Hence every recorded integer is at most \(L\), and
\(\ell(\sigma')\le L=\ell(\beta)\). Consequently, \(\nu\) has exactly \(L\) parts, all nonzero modulo \(k\).

The construction preserves weight, since
\[
|\nu|
=
|\beta|+k|\sigma'|
=
|\beta|+k|\sigma|
=
|\tau|
=
|\delta|.
\]
Together with the initial conjugation, this gives the bijection
\(\Theta_k:\mathcal D_k\to\mathcal O_k\).
Its inverse is described in \cite[Section~3]{xio19}.
In particular, \(\Theta_k(\varnothing)=\varnothing\).

\begin{exam}
\label{ex:stockhofe-keith}
\normalfont
Let \(k=3\) and \(\delta=(7,5,5,3,2)\). Then
\(\tau=\delta'=(5,5,4,3,3,1,1)\).
Step~1 deletes one occurrence of \(3\), recording \(1\).
The remaining \(3\) is the fourth part, so Step~2 deletes it
and subtracts \(3\) from the preceding three parts, recording
\(1+4-1=4\). Thus
\[
\begin{aligned}
(5,5,4,3,3,1,1)
&\xrightarrow{\text{Step 1}} (5,5,4,3,1,1)\\
&\xrightarrow{\text{Step 2}} (2,2,1,1,1).
\end{aligned}
\]
Therefore \(\beta=(2,2,1,1,1)\), \(\sigma=(4,1)\), and
\(\sigma'=(2,1,1,1)\), so
\[
\Theta_3(\delta)
=
(2,2,1,1,1)+3(2,1,1,1,0)
=
(8,5,4,4,1).
\]
Both \(\delta\) and \(\Theta_3(\delta)\) have weight \(22\).
\end{exam}

Read the parts of \(\tau=\delta'\) in weakly decreasing order,
reduce them modulo \(k\), and omit the zero residues.
Write the resulting sequence as \(v(\delta)=v_1v_2\cdots v_L\), where
\(1\le v_j\le k-1\).
Steps~1 and~2 delete only multiples of \(k\) and preserve
the residues and relative order of the surviving parts.
Step~3 adds multiples of \(k\) componentwise without
creating new parts. Hence
\begin{equation}
\label{eq:sec2-12}
\ell(\nu)=L,
\qquad
\nu_j\equiv v_j\pmod{k}
\quad (1\le j\le L).
\end{equation}
Thus the order of the nonzero residues is preserved,
not merely their multiplicities.

Combining the interpretation of \(\Sigma_i^{(k)}(\delta)\) as a residue-class count
in \(\tau\) with \eqref{eq:sec2-12} gives
\begin{equation}
\label{eq:sec2-14}
c_i^{(k)}(\nu)=\Sigma_i^{(k)}(\delta),
\qquad 1\le i<k,
\end{equation}
where \(c_i^{(k)}(\nu)\) counts the parts of \(\nu\)
congruent to \(i\pmod{k}\).
Moreover, each deletion removes one multiple of \(k\)
and contributes one part to \(\sigma\). Therefore \(\ell(\sigma)=C_k(\delta)=\Sigma_k^{(k)}(\delta)\).

Suppose that \(\delta\ne\varnothing\), and define
\(\operatorname{asc}(v)=\#\{j:1\le j<L,\ v_j<v_{j+1}\}\).
Since \(\beta\) is both \(k\)-flat and \(k\)-regular,
its residues determine it uniquely:
\(\beta_L=v_L\), and
\begin{equation}
\label{eq:sec2-17}
\beta_j-\beta_{j+1}
=
v_j-v_{j+1}
+
k\mathbf 1_{\{v_j<v_{j+1}\}},
\qquad 1\le j<L.
\end{equation}
Summing \eqref{eq:sec2-17} yields \(\beta_1=v_1+k\operatorname{asc}(v)\).

\begin{defi}
\label{def:residual-span}
\normalfont
For \(\delta\ne\varnothing\), define
\begin{equation*}
R_k(\delta)
=
k\bigl(C_k(\delta)+\operatorname{asc}(v(\delta))\bigr)+v_1,
\end{equation*}
and set \(R_k(\varnothing)=0\).
Using \(C_k(\delta)=\delta_1-s_k(\delta)\), we may equivalently write
\begin{equation}
\label{eq:sec2-20}
R_k(\delta)
=
k\bigl(\delta_1-s_k(\delta)
+\operatorname{asc}(v(\delta))\bigr)+v_1
\qquad(\delta\ne\varnothing).
\end{equation}
\end{defi}

Since \(\sigma'_1=\ell(\sigma)=C_k(\delta)\),
equation \eqref{eq:sec2-9} gives
\begin{equation}
\label{eq:sec2-21}
\Theta_k(\delta)_1=R_k(\delta),
\end{equation}
with the usual convention for the empty partition.
This is the largest-part refinement in \cite{fu19},
expressed in the present notation.

For \(k=2\) and \(\delta\ne\varnothing\), all nonzero residues
are \(1\), so \(\operatorname{asc}(v)=0\).
Writing \(s(\delta)=\Sigma_1^{(2)}(\delta)\), we obtain
\(R_2(\delta)=2\bigl(\delta_1-s(\delta)\bigr)+1\),
which agrees with the statistic used in the pairing-width
construction of Andrews and Dastidar~\cite{and26}.

\subsection{The bijection \texorpdfstring{\(\Phi_k\)}{Phi-k} and the statistic \texorpdfstring{\(W_k\)}{W-k}}
\label{subsec:global-correspondence}

We now return to an arbitrary partition \(\lambda\), with no restriction on its multiplicities.

For every positive integer \(j\), write uniquely \(m_j(\lambda) = ka_j+\varepsilon_j, \quad 0\leq\varepsilon_j<k\). Let \(\pi\) and \(\delta\) be the partitions defined by \(m_j(\pi)=a_j, \quad m_j(\delta)=\varepsilon_j\). We call \(\lambda\longleftrightarrow(\pi,\delta)\) the \(k\)-pairing decomposition of \(\lambda\). Equivalently,
\begin{equation}
\label{eq:sec2-24}
\lambda
=
\underbrace{
\pi\sqcup\pi\sqcup\cdots\sqcup\pi
}_{k\text{ copies}}
\sqcup\delta,
\end{equation}
and therefore \(|\lambda| = k|\pi|+|\delta|\). Each part of \(\pi\) records a group of \(k\)
equal parts of \(\lambda\), while \(m_j(\delta)=\varepsilon_j\) records the remaining copies of
\(j\).

Define \(r_k(\lambda)=\pi_1\), with \(r_k(\lambda)=0\) when \(\pi=\varnothing\). Thus
\(r_k(\lambda)\) is the largest integer \(j\) for which \(m_j(\lambda)\geq k\), or zero
if no such \(j\) exists.

We also define the \(k\)-pairing index and \(k\)-pairing rank by \(T_k(\lambda) = r_k(\lambda)+s_k(\delta)\) and \(\rho_k(\lambda) = r_k(\lambda)-s_k(\delta)\). Put
\[
P_k(\lambda)
=
k\ell(\pi).
\]

With \(R_k(\delta)=\Theta_k(\delta)_1\) as in \eqref{eq:sec2-21}, define the
\(k\)-pairing width by
\[
W_k(\lambda)
=
\max\{
P_k(\lambda),R_k(\delta)
\}.
\]

Thus, when \(\delta\neq\varnothing\),
\[
W_k(\lambda)
=
\max
\left\{
k\ell(\pi),
\,
k\bigl(
\delta_1-s_k(\delta)
+
\operatorname{asc}(v(\delta))
\bigr)+v_1
\right\}.
\]

If \(\delta=\varnothing\), the second entry in the maximum is defined to be \(0\).

For \(k=2\), \(P_2(\lambda) = 2\ell(\pi)\) and, for nonempty \(\delta\), \eqref{eq:sec2-20} gives \(R_2(\delta) =
2(\delta_1-s(\delta))+1\). Hence \(W_2(\lambda)=W(\lambda)\), so the definition recovers the pairing
width of Andrews and Dastidar~\cite{and26}.

\begin{defi}
\label{def:global-k-pairing}
\normalfont

Let \(\eta=\Theta_k(\delta)\). Define
\[
\Phi_k(\lambda)
=
k\pi'\sqcup\eta.
\]

Here \(k\pi'\) means that every part of \(\pi'\) is multiplied by \(k\), whereas \(k\pi'\sqcup\eta\)
means that the parts of \(k\pi'\) and \(\eta\) are combined and rearranged in weakly decreasing
order.

\end{defi}

The conjugation in \(k\pi'\) is essential. Indeed, \(\ell(\pi')=\pi_1, \quad \pi'_1=\ell(\pi)\), and therefore
\begin{equation}
\label{eq:sec2-33}
\ell(k\pi')=\pi_1=r_k(\lambda)
\end{equation}
and
\[
(k\pi')_1
=
k\ell(\pi)
=
P_k(\lambda).
\]

Thus \(k\pi'\) has \(\pi_1\) parts, each divisible by \(k\), and its largest part is
\(k\ell(\pi)\). In particular, the largest part of \(\pi\) becomes the number of parts of
\(k\pi'\), while \(k\) times the number of groups becomes its largest part.

\begin{thm}
\label{thm:k-pairing-width}
The map \(\Phi_k:\mathcal P\longrightarrow\mathcal P\) is a weight-preserving bijection.
Let \(\mu=\Phi_k(\lambda)\), and let \(c_i^{(k)}(\mu)\) count its parts congruent to
\(i\pmod k\), with multiplicity. Then \(c_0^{(k)}(\mu)=r_k(\lambda)\),
\(c_i^{(k)}(\mu)=\Sigma_i^{(k)}(\delta)\) for \(1\leq i<k\), and
\(\mu_1=W_k(\lambda)\). Hence, on partitions of each fixed size,
\[
\begin{aligned}
 &(r_k(\lambda),\Sigma_1^{(k)}(\delta),\ldots,
       \Sigma_{k-1}^{(k)}(\delta),W_k(\lambda))\\
 &\qquad\longleftrightarrow
 (c_0^{(k)}(\mu),c_1^{(k)}(\mu),\ldots,c_{k-1}^{(k)}(\mu),\mu_1).
\end{aligned}
\]
\end{thm}

\pf
We first prove bijectivity. Given \(\lambda\), the \(k\)-pairing decomposition
\(\lambda\longleftrightarrow(\pi,\delta)\) is unique by Euclidean division of each multiplicity.

Conversely, let \(\mu\) be any partition. Split it uniquely as \(\mu=\gamma\sqcup\eta\), where every
part of \(\gamma\) is divisible by \(k\), while \(\eta\in\mathcal O_k\). Then \(\pi=(\gamma/k)'\) is
uniquely determined, and since \(\Theta_k\) is a bijection, \(\delta=\Theta_k^{-1}(\eta)\) is
uniquely determined. The pair \((\pi,\delta)\) uniquely reconstructs \(\lambda\) by
\eqref{eq:sec2-24}. Hence \(\Phi_k\) is bijective.

For the weight,
\[
 |\Phi_k(\lambda)|=|k\pi'|+|\Theta_k(\delta)|=k|\pi|+|\delta|=|\lambda|.
\]
Thus \(\Phi_k\) is weight-preserving.

All parts of \(k\pi'\) are divisible by \(k\), whereas no part of \(\eta\) is divisible by \(k\). Hence the parts of \(\mu\) congruent to \(0\pmod k\) are precisely the parts of \(k\pi'\). By \eqref{eq:sec2-33}, \(c_0^{(k)}(\mu) = \ell(k\pi') = \pi_1 = r_k(\lambda)\). For \(1\leq i<k\), all parts of \(\mu\) congruent to \(i\pmod k\) come from \(\eta\). By \eqref{eq:sec2-14}, \(c_i^{(k)}(\mu) = c_i^{(k)}(\eta) = \Sigma_i^{(k)}(\delta)\). Finally,
\[
 \mu_1=\max\{(k\pi')_1,\eta_1\}=\max\{k\ell(\pi),R_k(\delta)\}=W_k(\lambda),
\]
by \eqref{eq:sec2-21} and the definition of \(W_k\). This proves the theorem. \qed

The theorem extends the pairing-width correspondence of Andrews and Dastidar~\cite{and26}.
Summing the residue-class counts gives \(\ell(\mu)=r_k(\lambda)+s_k(\delta)=T_k(\lambda)\),
and hence the following consequence.

\begin{core}
\label{cor:k-pairing-index-width}

For every partition \(\lambda\),
\[
\bigl(
T_k(\lambda),W_k(\lambda)
\bigr)
=
\bigl(
\ell(\Phi_k(\lambda)),
\Phi_k(\lambda)_1
\bigr).
\]

\end{core}

\subsection{Joint distributions and Gaussian polynomials}
\label{subsec:gaussian-enumeration}

Theorem~\ref{thm:k-pairing-width} gives the following generating functions. Let \(B=kh+r\), with \(h\geq0\) and \(0\leq r<k\). Since
\(W_k(\lambda)=\Phi_k(\lambda)_1\), the condition \(W_k(\lambda)\leq B\) corresponds to partitions
whose largest part is at most \(B\). Keeping track of the residue classes of the parts gives the
following.

\begin{thm}
\label{thm:bounded-width-product}

For \(B=kh+r\), \(0\leq r<k\),
\begin{equation}
\label{eq:sec2-40}
\begin{aligned}
&\sum_{\substack{\lambda\\ W_k(\lambda)\leq kh+r}}
y_0^{r_k(\lambda)}
\prod_{i=1}^{k-1}
y_i^{\Sigma_i^{(k)}(\delta)}
q^{|\lambda|}\\
&\qquad=
\frac{1}{(y_0q^k;q^k)_h}
\prod_{i=1}^{r}
\frac{1}{(y_iq^i;q^k)_{h+1}}
\prod_{i=r+1}^{k-1}
\frac{1}{(y_iq^i;q^k)_h}.
\end{aligned}
\end{equation}

\end{thm}

\pf
Under \(\Phi_k\), the exponent of \(y_i\) records the number of output parts congruent to \(i\pmod
k\). Among the positive integers not exceeding \(kh+r\), there are \(h\) multiples of \(k\), \(h+1\)
allowed parts in residue classes \(1,\ldots,r\), and \(h\) allowed parts in residue classes
\(r+1,\ldots,k-1\). Multiplying the corresponding geometric factors gives \eqref{eq:sec2-40}. \qed

Setting \(y_0=y_1=\cdots=y_{k-1}=z\) gives
\[
\sum_{\substack{\lambda\\W_k(\lambda)\leq B}}
z^{T_k(\lambda)}q^{|\lambda|}
=
\frac{1}{(zq;q)_B}.
\]

Letting \(B\to\infty\), we obtain
\[
\sum_{\lambda}
y^{r_k(\lambda)}
\prod_{i=1}^{k-1}
z_i^{\Sigma_i^{(k)}(\delta)}
q^{|\lambda|}
=
\frac{1}{(yq^k;q^k)_\infty}
\prod_{i=1}^{k-1}
\frac{1}{(z_iq^i;q^k)_\infty}.
\]

In particular, \(r_k\) and \(s_k\) are jointly equidistributed with the numbers of parts divisible and not divisible by \(k\), respectively. Setting all variables equal to \(z\) gives
\[
\sum_{\lambda}z^{T_k(\lambda)}q^{|\lambda|}
=
\frac{1}{(zq;q)_\infty}.
\]

Thus \(T_k\) is equidistributed with the number of parts.

We now refine simultaneously by \(T_k\) and \(W_k\). For \(t,b\geq1\), let \(N_{t,b}^{(k)}(n) = \#\{
\lambda\vdash n: T_k(\lambda)=t,\, W_k(\lambda)=b \}\).

\begin{core}
\label{cor:gaussian-pairing}

For \(t,b\geq1\),
\begin{equation}
\label{eq:sec2-47}
\sum_{n\geq0}
N_{t,b}^{(k)}(n)q^n
=
q^{t+b-1}
\begin{bmatrix}
t+b-2\\
t-1
\end{bmatrix}_q.
\end{equation}

\end{core}

\pf
By Corollary~\ref{cor:k-pairing-index-width}, \(N_{t,b}^{(k)}(n)\) counts partitions of \(n\) having
exactly \(t\) parts and largest part \(b\). Removing a first row of length \(b\) and one cell from
each of the remaining \(t-1\) rows leaves an arbitrary partition contained in a \((t-1)\times(b-1)\)
rectangle. Its generating function is the Gaussian polynomial in \eqref{eq:sec2-47}, while the
removed cells contribute \(q^{t+b-1}\). \qed

The right-hand side of \eqref{eq:sec2-47} is independent of \(k\). Hence the joint distribution of
\((T_k,W_k)\) is the same for every \(k\geq2\). For \(k=2\), \eqref{eq:sec2-47} recovers the
corresponding Gaussian-polynomial identity of Andrews and  Dastidar~\cite{and26}. Setting \(q=1\) gives
\(\sum_{n\geq0}N_{t,b}^{(k)}(n) = \binom{t+b-2}{t-1}\), so, for any prime \(p\), Kummer's theorem
identifies its \(p\)-adic valuation with the number of carries when \(t-1\) and \(b-1\) are added in
base \(p\), independently of \(k\).

\subsection{A bound on the largest part of the original partition}
\label{subsec:largest-part-bound}

In this subsection, the bound is imposed on \(\lambda_1\), rather than on \(W_k(\lambda)\).
Let \(\delta\in\mathcal D_k\) satisfy \(\delta_1\leq N\), and put \(\tau=\delta'\).
Pad \(\tau\) with zeros to length \(N\), and define \(w=w_1\cdots w_N\) by
\(w_j\equiv\tau_j\pmod k\) and \(0\leq w_j<k\). Unlike \(v(\delta)\), the sequence \(w\)
retains the zero residues, including those from the padding. Since \(\tau\) is \(k\)-flat,
\(w\) determines \(\tau\) uniquely:
\begin{equation}
\label{eq:sec2-51}
\tau_j
=
w_j+
k\#\{
a:j\leq a<N,\ w_a<w_{a+1}
\}.
\end{equation}

Indeed, flatness forces \(\tau_j-\tau_{j+1} = w_j-w_{j+1} + k\mathbf1_{\{w_j<w_{j+1}\}}\). Define \(\operatorname{amaj}(w)=\sum_{1\leq j<N,\,w_j<w_{j+1}}j\). Summing \eqref{eq:sec2-51} gives
\begin{equation}
\label{eq:sec2-54}
|\delta|
=
\sum_{j=1}^{N}w_j
+
k\operatorname{amaj}(w).
\end{equation}

If the letter \(i\) occurs \(n_i\) times in \(w\), then \(\Sigma_i^{(k)}(\delta)=n_i \quad(1\leq i<k)\). We use the following \(q\)-multinomial identity for multiset words~\cite{sta12}:
\begin{equation}
\label{eq:sec2-57}
\sum_{\substack{
w\in\{0,\ldots,k-1\}^N\\
\#\{j:w_j=i\}=n_i
}}
q^{k\operatorname{amaj}(w)}
=
\begin{bmatrix}
N\\
n_0,n_1,\ldots,n_{k-1}
\end{bmatrix}_{q^k},
\end{equation}
where \(n_0+\cdots+n_{k-1}=N\). 

Combining \eqref{eq:sec2-54} and \eqref{eq:sec2-57} gives the following generating function
for partitions \(\delta\in\mathcal D_k\) with \(\delta_1\leq N\):
\begin{equation}
\label{eq:sec2-62}
\begin{aligned}
&
\sum_{\substack{
\delta\in\mathcal D_k\\
\delta_1\leq N
}}
\prod_{i=1}^{k-1}
z_i^{\Sigma_i^{(k)}(\delta)}
q^{|\delta|}
=
\sum_{\substack{
t_1,\ldots,t_{k-1}\geq0\\
T=t_1+\cdots+t_{k-1}\leq N
}}
\begin{bmatrix}
N\\
N-T,t_1,\ldots,t_{k-1}
\end{bmatrix}_{q^k}
\prod_{i=1}^{k-1}
(z_iq^i)^{t_i}.
\end{aligned}
\end{equation}

For an arbitrary partition \(\lambda\longleftrightarrow(\pi,\delta)\), the condition \(\lambda_1\leq N\) is equivalent to \(\pi_1\leq N,\quad \delta_1\leq N\). For fixed \(\pi_1=r\geq0\), with the empty partition corresponding to \(r=0\), \(\sum_{\pi:\pi_1=r}q^{k|\pi|} = \frac{q^{kr}}{(q^k;q^k)_r}\). Hence \eqref{eq:sec2-62} immediately yields
\begin{equation}
\label{eq:sec2-63}
\begin{aligned}
&
\sum_{\substack{\lambda\\\lambda_1\leq N}}
y^{r_k(\lambda)}
\prod_{i=1}^{k-1}
z_i^{\Sigma_i^{(k)}(\delta)}
q^{|\lambda|}\\
&=
\left(
\sum_{r=0}^{N}
\frac{y^rq^{kr}}
{(q^k;q^k)_r}
\right)
\sum_{\substack{
t_1,\ldots,t_{k-1}\geq0\\
T\leq N
}}
\begin{bmatrix}
N\\
N-T,t_1,\ldots,t_{k-1}
\end{bmatrix}_{q^k}
\prod_{i=1}^{k-1}
(z_iq^i)^{t_i}.
\end{aligned}
\end{equation}

To combine the nonzero residue statistics, for \(t\geq0\) define
\[
H_{k,t}(q)
=
\sum_{\substack{
t_1,\ldots,t_{k-1}\geq0\\
t_1+\cdots+t_{k-1}=t
}}
q^{\sum_{i=1}^{k-1}it_i}
\begin{bmatrix}
t\\
t_1,\ldots,t_{k-1}
\end{bmatrix}_{q^k}.
\]

Then, with
\[
P_N^{(k)}(y,z,q)
=
\sum_{\substack{\lambda\\\lambda_1\leq N}}
y^{r_k(\lambda)}
z^{s_k(\delta)}
q^{|\lambda|},
\]
equation \eqref{eq:sec2-63} becomes
\begin{equation}
\label{eq:sec2-66}
P_N^{(k)}(y,z,q)
=
\sum_{r=0}^{N}\sum_{t=0}^{N}
\frac{
y^rz^tq^{kr}H_{k,t}(q)
}{
(q^k;q^k)_r
}
\begin{bmatrix}
N\\t
\end{bmatrix}_{q^k}.
\end{equation}

For \(k=2\), \(H_{2,t}(q)=q^t\), so \eqref{eq:sec2-66} reduces to the finite generating function in
the original pairing theory~\cite[Lemma~2.2]{and26}.

Finally, for \(t\geq0\), define
\[
B_{k,t}(q)
=
\sum_{\substack{
\delta\in\mathcal D_k\\
s_k(\delta)=t
}}
q^{|\delta|}.
\]

Letting \(N\to\infty\) in \eqref{eq:sec2-62} and extracting total degree \(t\) gives
\begin{equation}
\label{eq:sec2-67}
B_{k,t}(q)
=
\sum_{\substack{
t_1,\ldots,t_{k-1}\geq0\\
t_1+\cdots+t_{k-1}=t
}}
\frac{
q^{\sum_{i=1}^{k-1}it_i}
}{
\prod_{i=1}^{k-1}
(q^k;q^k)_{t_i}
}
=
\frac{H_{k,t}(q)}
{(q^k;q^k)_t}.
\end{equation}

This form will be used in Section~\ref{sec:negative-rank}.

\section{Negative \texorpdfstring{$k$}{k}-Pairing Rank and Marked Gaps}
\label{sec:negative-rank}

In the original pairing theory, Andrews and Dastidar~\cite{and26} call a partition simply paired
if no part occurs more than three times. For such a partition, every part size contributes at most
one pair. Hence the partition \(\pi\) recording the pairs has distinct parts and may be
identified with the finite set of its part sizes. They use this observation to
study simply paired partitions of negative pairing rank and obtain signed and unsigned enumerations
involving odd divisors and overpartitions into odd parts.

We now develop the corresponding theory for \(k\)-pairing.

\subsection{The bijection \texorpdfstring{\(\Psi_k\)}{Psi-k}}
\label{subsec:simply-k-paired}
\label{subsec:gap-correspondence}

\begin{defi}
\label{def:simply-k-paired}
\normalfont
A partition \(\lambda\) is \emph{simply \(k\)-paired} if \(m_j(\lambda)\leq2k-1\) for every
\(j\geq1\). Let \(\mathcal S_k\) denote the set of these partitions, and let
\(\mathcal S_k^-:=\{\lambda\in\mathcal S_k:\rho_k(\lambda)<0\}\).
\end{defi}
In the \(k\)-pairing decomposition, the condition \(\lambda\in\mathcal S_k\) is equivalent
to \(a_j\in\{0,1\}\) for every \(j\). The groups of \(k\) equal parts are therefore indexed by
\(A=\{j:m_j(\lambda)\geq k\}\). Thus \(\lambda\longleftrightarrow(A,\delta)\) is unique, with
\(m_j(\lambda)=k\mathbf1_{\{j\in A\}}+m_j(\delta)\).
When taking a conjugate, we regard \(A\) as the partition obtained by listing its elements
in decreasing order. Put \(b_k(\lambda)=|A|\), the number of groups of \(k\) equal parts,
and write \(s_k(\lambda)=s_k(\delta)\). Then
\[
  |\lambda|=k\sum_{a\in A}a+|\delta|,\qquad
  r_k(\lambda)=\max A,\qquad \rho_k(\lambda)=\max A-s_k(\delta).
\]
With \(\max\varnothing=0\) and \(t=s_k(\delta)\),
\begin{equation}
\label{eq:sec3-10}
  \lambda\in\mathcal S_k^-\quad\Longleftrightarrow\quad
  t\geq1,\qquad A\subseteq\{1,\ldots,t-1\}.
\end{equation}

For a nonempty partition \(\nu=(\nu_1,\ldots,\nu_t)\), define
\[
G_k(\nu)=\{j:1\leq j<t,\ \nu_j-\nu_{j+1}\geq k\},\qquad
  g_k(\nu)=|G_k(\nu)|.
\]
Only internal gaps are counted, not the gap below the smallest positive part. Set
\(\eta=\Theta_k(\delta)\). Since \(\ell(\eta)=s_k(\delta)=t\) and
\(A'_i-A'_{i+1}=\mathbf1_{\{i\in A\}}\), the componentwise sum \(\nu=\eta+kA'\) satisfies
\begin{equation}
\label{eq:sec3-27}
  \nu_i-\nu_{i+1}=\eta_i-\eta_{i+1}+k\mathbf1_{\{i\in A\}}.
\end{equation}
A group of \(k\) parts equal to \(a\) adds \(k\) to each of the first \(a\) parts of
\(\eta\), increasing only the \(a\)-th internal gap. Unlike the map \(\Phi_k\), the map
\(\Psi_k\) records the selected positions as part of its output:
\[
\Psi_k(\lambda)=(\nu,J):=(\eta+kA',A).
\]

\begin{thm}
\label{thm:marked-gap-bijection}
The map \(\Psi_k\) is a weight-preserving bijection
\[
\mathcal S_k^-\longleftrightarrow
  \{(\nu,J):\nu\in\mathcal O_k\setminus\{\varnothing\},\ J\subseteq G_k(\nu)\}.
\]
If \(\Psi_k(\lambda)=(\nu,J)\), then
\[
s_k(\lambda)=\ell(\nu),\qquad b_k(\lambda)=|J|,\qquad r_k(\lambda)=\max J.
\]
\end{thm}
\pf
Let \(\lambda\longleftrightarrow(A,\delta)\) be a simply \(k\)-paired partition of negative
\(k\)-pairing rank. Put \(\eta=\Theta_k(\delta), \quad t=\ell(\eta)=s_k(\delta)\). By
\eqref{eq:sec3-10}, \(A\subseteq\{1,\ldots,t-1\}\). Define \(\nu=\eta+kA', \quad J=A\). For \(1\leq
i<t\), equation \eqref{eq:sec3-27} gives \(\nu_i-\nu_{i+1} = \eta_i-\eta_{i+1} + k\mathbf1_{\{i\in
A\}} \geq0\). Moreover, since \(A\subseteq\{1,\ldots,t-1\}\), \(A'_t=0\), and therefore
\(\nu_t=\eta_t>0\). Thus \(\nu\) has exactly \(t\) positive parts.

Each part of \(\nu\) differs from the corresponding part of \(\eta\) by a multiple of \(k\). Hence \(\nu_i\equiv\eta_i\pmod k\). Since \(\eta\in\mathcal O_k\), \(k\nmid\eta_i\), and therefore also \(k\nmid\nu_i\). Thus \(\nu\in\mathcal O_k\). If \(i\in A\), then \eqref{eq:sec3-27} yields \(\nu_i-\nu_{i+1}\geq k\), so \(J=A\subseteq G_k(\nu)\). The map is weight-preserving because
\[
|\nu|=|\eta|+k|A'|=|\delta|+k\sum_{a\in A}a=|\lambda|.
\]

Here we used \(|A'|=\sum_{a\in A}a\). It remains to construct the
inverse.

Let \((\nu,J)\) satisfy \(\nu\in\mathcal O_k\setminus\{\varnothing\}, \quad J\subseteq G_k(\nu)\). Write \(t=\ell(\nu)\) and define
\[
\eta=\nu-kJ'.
\]

For \(1\leq i<t\),
\begin{equation}
\label{eq:sec3-36}
\eta_i-\eta_{i+1}
=
\nu_i-\nu_{i+1}
-
k\mathbf1_{\{i\in J\}}.
\end{equation}

If \(i\notin J\), the right-hand side is nonnegative because \(\nu\) is a partition.

If \(i\in J\), then \(i\in G_k(\nu)\), so \(\nu_i-\nu_{i+1}\geq k\), and \eqref{eq:sec3-36} is again
nonnegative.

Furthermore, \(J\subseteq\{1,\ldots,t-1\}\), so \(J'_t=0\) and hence \(\eta_t=\nu_t>0\). Thus
\(\eta\) is a partition with exactly \(t\) positive parts.

Since only multiples of \(k\) have been subtracted, \(\eta_i\equiv\nu_i\pmod k\). Because
\(\nu\in\mathcal O_k\), we also have \(\eta\in\mathcal O_k\). The bijection \(\Theta_k\) therefore
uniquely recovers \(\delta=\Theta_k^{-1}(\eta)\). Set \(A=J\). Then \(\max A<t=s_k(\delta)\), so the
simply \(k\)-paired partition determined by \((A,\delta)\) has negative \(k\)-pairing rank.

The two constructions are plainly inverse to one another. Finally, \(s_k(\lambda) = s_k(\delta) =
\ell(\eta) = \ell(\nu)\), \(b_k(\lambda) = |A| = |J|\), and \(r_k(\lambda) = \max A = \max J\). This
proves the theorem. \qed

For example, \(k=3\), \(A=(3,1)\), and \(\eta=(2,2,2,2)\) give
\(\nu=(2,2,2,2)+3(2,1,1,0)=(8,5,5,2)\), with \(J=\{1,3\}\).

\begin{core}
\label{cor:marked-gap-series}
We have
\begin{equation}
\label{eq:sec3-17}
\begin{aligned}
 F_k(u;q)&:=\sum_{\lambda\in\mathcal S_k^-}u^{b_k(\lambda)}q^{|\lambda|}
 =\sum_{\nu\in\mathcal O_k\setminus\{\varnothing\}}(1+u)^{g_k(\nu)}q^{|\nu|}\\
 &=\sum_{t\geq1}B_{k,t}(q)(-uq^k;q^k)_{t-1}
 =\sum_{t\geq1}\frac{H_{k,t}(q)(-uq^k;q^k)_{t-1}}{(q^k;q^k)_t}.
\end{aligned}
\end{equation}
\end{core}
\pf
For fixed \(\nu\), summing over \(J\subseteq G_k(\nu)\) gives
\((1+u)^{g_k(\nu)}\). For fixed \(s_k(\delta)=t\), the independent subsets in
\eqref{eq:sec3-10} contribute
\(\prod_{j=1}^{t-1}(1+uq^{kj})=(-uq^k;q^k)_{t-1}\).
The last equality is \eqref{eq:sec2-67}. \qed

\subsection{Signed and unsigned enumerations}
\label{subsec:gap-enumerations}

Write \(\mathcal O_k(n)=\{\nu\in\mathcal O_k:|\nu|=n\}\). Let \(E_k^{\mathrm{sp}}(n)\) and \(O_k^{\mathrm{sp}}(n)\) count partitions in
\(\mathcal S_k^-\) of size \(n\) with \(b_k\) even and odd, respectively, and put
\(M_k(n)=E_k^{\mathrm{sp}}(n)+O_k^{\mathrm{sp}}(n)\).

\begin{core}
\label{thm:signed-gap-cancellation}
\label{thm:unsigned-gap-enumeration}
For \(n\geq1\),
\[
\begin{aligned}
 E_k^{\mathrm{sp}}(n)-O_k^{\mathrm{sp}}(n)
 &=\#\{\nu\in\mathcal O_k(n):g_k(\nu)=0\},\\
 M_k(n)&=\sum_{\nu\in\mathcal O_k(n)}2^{g_k(\nu)}.
\end{aligned}
\]
Their generating functions are
\[
\begin{aligned}
 \sum_{n\geq1}\bigl(E_k^{\mathrm{sp}}(n)-O_k^{\mathrm{sp}}(n)\bigr)q^n
 &=\sum_{t\geq1}\frac{H_{k,t}(q)}{1-q^{kt}},\\
 \sum_{n\geq1}M_k(n)q^n
 &=\sum_{t\geq1}\frac{H_{k,t}(q)(-q^k;q^k)_{t-1}}{(q^k;q^k)_t}.
\end{aligned}
\]
\end{core}
\pf
Set \(u=-1\) and \(u=1\) in \eqref{eq:sec3-17}. In the signed case, fix \(\nu\).
If \(G_k(\nu)\neq\varnothing\), let \(j=\min G_k(\nu)\). Replacing \(J\) by
\(J\mathbin{\triangle}\{j\}\) reverses the sign and pairs all subsets of \(G_k(\nu)\).
Only the terms with \(g_k(\nu)=0\) remain. In the unsigned case there are
\(2^{g_k(\nu)}\) subsets. The series follow using
\((q^k;q^k)_{t-1}/(q^k;q^k)_t=(1-q^{kt})^{-1}\). \qed

The surviving partitions have \(k\nmid\nu_j\) and
\(\nu_j-\nu_{j+1}<k\) for \(1\leq j<\ell(\nu)\). They need not be \(k\)-flat, since their
smallest part is unrestricted. The divisor and overpartition interpretations for \(k=2\)
are included in the following results, obtained by requiring every part of \(\Theta_k(\delta)\)
to lie in one fixed nonzero residue class modulo \(k\).

\subsection{Restriction to one nonzero residue class}
\label{subsec:divisors-overpartitions}

Fix \(1\leq a<k\). Let \(\mathcal S_{k,a}^-\) consist of the partitions
\(\lambda\in\mathcal S_k^-\) for which the partition \(\delta\) in the decomposition
\(\lambda\longleftrightarrow(A,\delta)\) satisfies \(\Sigma_i^{(k)}(\delta)=0\) for
\(1\leq i<k\), \(i\neq a\). Equivalently, every part of \(\eta=\Theta_k(\delta)\) is
congruent to \(a\pmod k\). The same is true of \(\nu=\eta+kA'\).
Recall that an overpartition is a partition in which the first occurrence of each distinct part
may be overlined~\cite{cor04}.
Use \(E_{k,a}^{\mathrm{sp}}\), \(O_{k,a}^{\mathrm{sp}}\), and \(M_{k,a}\) for the corresponding
restricted counts, and let \(\overline p_{a,k}(n)\) count overpartitions with all parts in this
residue class.

If \(\eta=\Theta_k(\delta)\) has \(t\) parts, write \(\eta_j=a+k\alpha_j\), where
\(\alpha_1\geq\cdots\geq\alpha_t\geq0\). Since \(|\eta|=at+k|\alpha|\), its generating
function is
\begin{equation}
\label{eq:sec3-48}
 B_{k,a,t}(q)
 :=\sum_{\substack{\delta\in\mathcal D_k,\ s_k(\delta)=t\\
 \Sigma_i^{(k)}(\delta)=0\ (1\leq i<k,\ i\neq a)}}q^{|\delta|}
 =\frac{q^{at}}{(q^k;q^k)_t}.
\end{equation}

\begin{thm}
\label{thm:residue-class-divisors}
For \(n\geq1\),
\begin{equation}
\label{eq:sec3-49}
 E_{k,a}^{\mathrm{sp}}(n)-O_{k,a}^{\mathrm{sp}}(n)
 =\#\{d\mid n:d\equiv a\pmod k\}.
\end{equation}
\end{thm}
\pf
By Corollary~\ref{thm:signed-gap-cancellation}, after signed cancellation only those \(k\)-regular
partitions \(\nu\) with \(g_k(\nu)=0\) remain.

Every part of \(\nu\) is congruent to \(a\pmod k\). Hence two distinct part sizes differ by a
positive multiple of \(k\), and therefore by at least \(k\).

But \(g_k(\nu)=0\) requires every internal gap to be strictly less than \(k\). It follows that no
two adjacent parts can be distinct. Hence \(\nu=(d^t)\) for some \(d\equiv a\pmod k\). Its weight is
\(n=dt\). Thus \(d\) is a divisor of \(n\) satisfying \(d\equiv a\pmod k\). Conversely, every such
divisor gives the unique rectangular partition \((d^{\,n/d})\). This proves \eqref{eq:sec3-49}. \qed

The equivalent series identity is
\[
\sum_{n\geq1}\bigl(E_{k,a}^{\mathrm{sp}}(n)-O_{k,a}^{\mathrm{sp}}(n)\bigr)q^n
 =\sum_{t\geq1}\frac{q^{at}}{1-q^{kt}}
 =\sum_{t\geq1}\sum_{j\geq0}q^{(a+kj)t}.
\]

\Needspace{11\baselineskip}
\begin{thm}
\label{thm:residue-class-overpartitions}
For \(n\geq1\),
\begin{equation}
\label{eq:sec3-51}
 M_{k,a}(n)=\frac12\,\overline p_{a,k}(n),
\end{equation}
and
\[
\sum_{n\geq1}M_{k,a}(n)q^n
 =\frac12\left(\frac{(-q^a;q^k)_\infty}{(q^a;q^k)_\infty}-1\right).
\]
\end{thm}
\pf
Let \(\nu\) be a nonempty partition all of whose parts are congruent to \(a\pmod k\), and let
\(d(\nu)\) denote its number of distinct part sizes.

Whenever the part size changes, the difference is a positive multiple of \(k\), and hence is at
least \(k\). Therefore \(g_k(\nu) = d(\nu)-1\). By the bijection \(\Psi_k\), there are
\(2^{d(\nu)-1}\) partitions \(\lambda\in\mathcal S_{k,a}^-\) whose images have first component
\(\nu\).

On the other hand, an overpartition with underlying partition \(\nu\) has two choices for each of
its \(d(\nu)\) distinct part sizes, giving \(2^{d(\nu)}\) overpartitions.

More explicitly, a marked gap after a block of equal parts may be identified with the choice to
overline that part size. Every distinct part size except the smallest has such an internal gap below
its block. Thus the pairs \((\nu,J)\) with \(J\subseteq G_k(\nu)\) correspond exactly to
those overpartitions in which the smallest part is not overlined.

Consequently, \(2^{g_k(\nu)} = 2^{d(\nu)-1} = \frac12\,2^{d(\nu)}\). Summing over all nonempty
\(\nu\) proves \eqref{eq:sec3-51}. The product identity follows by assigning the factor
\((1+q^{a+kj})/(1-q^{a+kj})\) to each allowed part size, multiplying these factors, removing the empty object, and taking
one half. \qed

The series form from \eqref{eq:sec3-17} is
\(\sum_{t\geq1}q^{at}(-q^k;q^k)_{t-1}/(q^k;q^k)_t\); its product evaluation also follows from the
\(q\)-binomial theorem. At \((k,a)=(2,1)\), the two theorems give
\(E_2^{\mathrm{sp}}(n)-O_2^{\mathrm{sp}}(n)=d_{\mathrm{odd}}(n)\) and
\(M_2(n)=\overline p_{\mathrm{odd}}(n)/2\), recovering~\cite[Section~6]{and26}.

\subsection{Generating functions for fixed \texorpdfstring{\(k\)}{k}-pairing rank}
\label{subsec:fixed-rank}

We conclude this section by fixing \(\rho_k(\lambda)=-d\), with \(d\geq0\).
Andrews and Dastidar~\cite{and26} introduced generating functions for fixed pairing rank
and derived a recurrence relating the cases \(\rho_2=-d\) and \(\rho_2=-(d+1)\).
For general \(k\), we obtain a corresponding recurrence under the condition that every part
of \(\Theta_k(\delta)\) is congruent to \(a\pmod k\), where \(1\leq a<k\).

For \(d\geq0\), let
\[
A_d^{(k)}(q)=\sum_{\substack{\lambda\in\mathcal S_k\\\rho_k(\lambda)=-d}}q^{|\lambda|}.
\]
Write \(A_d^{(k,a)}(q)\) for the same generating function with the additional requirement
that every part of \(\Theta_k(\delta)\) is congruent to \(a\pmod k\).
The case \(\rho_k=0\) includes the empty partition.

\begin{pro}
\label{prop:fixed-rank-series}
For \(d\geq0\),
\begin{equation}
\label{eq:sec3-57}
 A_d^{(k)}(q)=B_{k,d}(q)
 +\sum_{r\geq1}q^{kr}(-q^k;q^k)_{r-1}B_{k,r+d}(q).
\end{equation}
In particular,
\begin{equation}
\label{eq:sec3-58}
 A_d^{(k,a)}(q)=\frac{q^{ad}}{(q^k;q^k)_d}
 +\sum_{r\geq1}\frac{q^{(k+a)r+ad}(-q^k;q^k)_{r-1}}{(q^k;q^k)_{r+d}}.
\end{equation}
\end{pro}
\pf
Fix \(r=\max A\). If \(r=0\), then \(A=\varnothing\) and \(s_k(\delta)=d\).
If \(r\geq1\), the set \(A\) contains \(r\) and an arbitrary subset of
\(\{1,\ldots,r-1\}\), contributing \(q^{kr}(-q^k;q^k)_{r-1}\), while
\(s_k(\delta)=r+d\). Summing gives \eqref{eq:sec3-57}, and
\eqref{eq:sec3-48} gives \eqref{eq:sec3-58}. \qed

\begin{thm}
\label{thm:fixed-rank-recurrence}
Fix \(k\geq2\) and \(1\leq a<k\). For every \(d\geq0\),
\begin{equation}
\label{eq:sec3-59}
\begin{aligned}
 (1+q^{k(d+1)})A_{d+1}^{(k,a)}(q)
 &=(q^a+q^{kd})A_d^{(k,a)}(q)
 -\frac{q^{(k+a)d}}{(q^k;q^k)_d}
 +\frac{q^{(k+a)(d+1)}}{(q^k;q^k)_{d+1}}.
\end{aligned}
\end{equation}
\end{thm}
\pf
Suppress the superscripts and put \(a_0=1\), \(a_r=q^{kr}(-q^k;q^k)_{r-1}\) for \(r\geq1\), and
\(b_t=q^{at}/(q^k;q^k)_t\). Then \eqref{eq:sec3-58} reads
\(A_d=\sum_{r\geq0}a_rb_{r+d}\). Since \((1-q^{k(t+1)})b_{t+1}=q^ab_t\),
\[
 A_{d+1}-q^{k(d+1)}S_d=q^aA_d,\qquad
 S_d:=\sum_{r\geq0}q^{kr}a_rb_{r+d+1}.
\]
For \(r\geq1\), \(q^{kr}a_r=a_{r+1}/q^k-a_r\). Separating \(r=0\) and using
\(a_0=1\), \(a_1=q^k\), we obtain
\[
 S_d=b_{d+1}+\frac{A_d-b_d-q^kb_{d+1}}{q^k}-(A_{d+1}-b_{d+1})
 =\frac{A_d-b_d}{q^k}-A_{d+1}+b_{d+1}.
\]
Substitution gives
\((1+q^{k(d+1)})A_{d+1}=(q^a+q^{kd})A_d-q^{kd}b_d+q^{k(d+1)}b_{d+1}\), which is
\eqref{eq:sec3-59}. \qed

For \(k=2\) and \(a=1\), this is the recurrence in~\cite[Lemma~7.1]{and26}.

\section{Equivalence classes of tuples of partitions}
\label{sec:partition-tuples}

The diagonal pairing of Andrews and Dastidar~\cite{and26} superimposes the two wings of a
Durfee square, with the lower wing reflected across the diagonal. Their common cells are
paired, and the edge-connected components of the symmetric difference may be transferred
independently between the two diagrams. We extend this construction to tuples of \(k\) Young
diagrams.

\subsection{An equivalence relation on tuples of partitions}
\label{subsec:tuple-equivalence}
\label{subsec:class-representatives}

\begin{defi}
\label{def:partition-tuples}
\normalfont
For \(k\geq2\) and \(d\geq0\), let \(\mathscr W_{k,d}\) be the set of ordered tuples
\[
 \boldsymbol\alpha=(d;\alpha^{(0)},\ldots,\alpha^{(k-1)}),\qquad
 \ell(\alpha^{(s)})\leq d\quad(0\leq s<k),
\]
where each \(\alpha^{(s)}\) is a partition. Put
\(\mathscr W_k=\bigsqcup_{d\geq0}\mathscr W_{k,d}\), and define the size of a tuple by
\[
N(\boldsymbol\alpha)=d^2+\sum_{s=0}^{k-1}|\alpha^{(s)}|.
\]
\end{defi}
For \(k=2\), the two component partitions reconstruct an ordinary partition with Durfee
size \(d\). For \(k>2\), we work with the tuples themselves, not with a decomposition of an
ordinary Young diagram. We identify a Young diagram with an order ideal in the product order
on \(\mathbb Z_{>0}^2\): \((i,j)\preceq(i',j')\) if \(i\leq i'\) and \(j\leq j'\).

For \(\boldsymbol\alpha\in\mathscr W_{k,d}\), define the cell multiplicity
\(h_{\boldsymbol\alpha}\) and the diagrams \(U_r\) by
\[
h_{\boldsymbol\alpha}(x)=\#\{s:x\in\alpha^{(s)}\},\qquad
 U_r=\{x:h_{\boldsymbol\alpha}(x)\geq r\}\quad(1\leq r\leq k).
\]
Each component diagram is an order ideal, so each \(U_r\) is also a Young diagram.
Moreover, \(U_1\supseteq\cdots\supseteq U_k\), with \(\ell(U_1)\leq d\).
For \(k=2\), the diagrams \(U_1,U_2\) are the union and intersection of the two components.

We use~\cite[Lemma~4.2]{and26}. For two Young diagrams \(B,C\), each edge-connected component
of \((B\cup C)\setminus(B\cap C)\) belongs to exactly one of the two diagrams and may be
transferred to the other without changing the union or intersection. On \(\mathscr W_{k,d}\),
let \(\sim\) be the equivalence relation generated by applying such transfers to any two
component diagrams while leaving the others unchanged. In particular, finitely many transfers
give the transformation
\[
(B,C)\longmapsto(B\cup C,B\cap C).
\]

\begin{thm}
\label{thm:unique-representative}
For fixed \(d\) and \(\boldsymbol\alpha,\boldsymbol\beta\in\mathscr W_{k,d}\),
\[
\boldsymbol\alpha\sim\boldsymbol\beta
 \quad\Longleftrightarrow\quad h_{\boldsymbol\alpha}=h_{\boldsymbol\beta}.
\]
Every equivalence class contains exactly one representative whose component diagrams
satisfy \(\alpha^{(0)}\supseteq\cdots\supseteq\alpha^{(k-1)}\), namely
\((d;U_1,\ldots,U_k)\).
\end{thm}
\pf
The function \(h_{\boldsymbol\alpha}\) is invariant under the transfers defining \(\sim\),
which proves the forward implication.

For the converse, we first show that every tuple in \(\mathscr W_{k,d}\) can be transformed
into one whose component diagrams are weakly decreasing under inclusion.
Consider two adjacent components \(\alpha^{(s)},\alpha^{(s+1)}\). By the lemma recalled above,
a sequence of transfers replaces them by
\begin{equation}
\label{eq:sec4-adjacent-normalization}
\left(
\alpha^{(s)}\cup\alpha^{(s+1)},
\,
\alpha^{(s)}\cap\alpha^{(s+1)}
\right).
\end{equation}

Define \(I(\boldsymbol{\alpha}) = \sum_{s=0}^{k-1} s\,|\alpha^{(s)}|\). Under \eqref{eq:sec4-adjacent-normalization}, the contribution of the two components changes from \(s|\alpha^{(s)}| + (s+1)|\alpha^{(s+1)}|\) to \(s|\alpha^{(s)}\cup\alpha^{(s+1)}| + (s+1)|\alpha^{(s)}\cap\alpha^{(s+1)}|\). Using \(|A\cup B|+|A\cap B|=|A|+|B|\), the decrease in \(I\) is \(\left| \alpha^{(s+1)} \setminus \alpha^{(s)} \right|\). Thus \(I\) decreases strictly whenever \(\alpha^{(s+1)} \not\subseteq \alpha^{(s)}\). Apply \eqref{eq:sec4-adjacent-normalization} whenever
\(\alpha^{(s+1)}\not\subseteq\alpha^{(s)}\). Each step strictly decreases the nonnegative
integer \(I\), so the procedure terminates. At termination,
\begin{equation}
\label{eq:sec4-9}
\alpha^{(0)}
\supseteq
\alpha^{(1)}
\supseteq
\cdots
\supseteq
\alpha^{(k-1)}.
\end{equation}

Now fix a cell \(x\) with \(h_{\boldsymbol\alpha}(x)=r\). In a tuple satisfying
\eqref{eq:sec4-9}, the component diagrams containing \(x\) are precisely the first \(r\)
components. Hence the component with index \(s\) is \(U_{s+1}\), and the terminal tuple is
\((d;U_1,\ldots,U_k)\). This also proves uniqueness.

If two tuples have the same multiplicity function \(h\), they have the same diagrams \(U_r\).
Both can therefore be transformed into the representative \((d;U_1,\ldots,U_k)\).
Since every transfer is reversible, the two tuples belong to the same equivalence class. \qed

For \(k=2\), the representative \((U_1,U_2)=(U,I)\) corresponds to the partition whose
successive ranks are all nonnegative, as in~\cite{and26}.

\subsection{Cardinalities of equivalence classes}
\label{subsec:class-cardinalities}

Fix a multiplicity function \(h=h_{\boldsymbol\alpha}\) arising from a tuple in
\(\mathscr W_{k,d}\), and put \(U_r=\{x:h(x)\geq r\}\).
Let \(\mathcal B_r\) be the set of edge-connected components of
\(L_r=U_r\setminus U_{r+1}=\{x:h(x)=r\}\), for \(1\leq r<k\), and put
\(\mathcal B=\bigsqcup_{r=1}^{k-1}\mathcal B_r\). Write \(r(C)=r\) for \(C\in\mathcal B_r\).
For any tuple \(\boldsymbol\alpha\in\mathscr W_{k,d}\) satisfying \(h_{\boldsymbol\alpha}=h\),
the subsets \(S(x)=\{s:x\in\alpha^{(s)}\}\) satisfy
\begin{equation}
\label{eq:sec4-14}
 |S(x)|=h(x),\qquad x\preceq y\ \Longrightarrow\ S(y)\subseteq S(x).
\end{equation}
If \(x,y\) are adjacent cells with \(h(x)=h(y)\), then one of \(S(x),S(y)\) contains
the other. Since they have the same cardinality, they are equal. Thus \(S\) is constant on
each component \(C\in\mathcal B\).

For \(C,D\in\mathcal B\), write \(C\succeq D\) if there exist \(x\in C\) and \(y\in D\)
with \(x\preceq y\), and take the reflexive transitive closure. This is a partial order:
\(r(C)\geq r(D)\) whenever \(C\succeq D\). Moreover, comparable cells on which \(h\) has the
same value can be joined by a path of horizontal and vertical steps on which \(h\) remains
constant, so they lie in the same edge-connected component. Define
\begin{equation}
\label{eq:sec4-15}
 \Omega_k(U_1,\ldots,U_k)
 =\#\{S:\mathcal B\to2^{\{0,\ldots,k-1\}}:
 |S(C)|=r(C),\ C\succeq D\Rightarrow S(D)\subseteq S(C)\}.
\end{equation}
For \(x\in U_k\), put \(S(x)=\{0,\ldots,k-1\}\); for \(x\notin U_1\), put
\(S(x)=\varnothing\).
We take \(b=0\) when \(U_1=\varnothing\), with the empty product equal to \(1\).

\begin{thm}
\label{thm:class-cardinality}
\label{prop:local-class-product}
Let \(\mathcal C_{d,h}\) be the equivalence class under \(\sim\) with multiplicity function
\(h\), set
\(b=(U_1)_1\), and write \(h_{i,j}=h(i,j)\). With the boundary values
\(h_{i,0}=k\) and \(h_{d+1,j}=0\),
\begin{equation}
\label{eq:sec4-20}
 |\mathcal C_{d,h}|=\Omega_k(U_1,\ldots,U_k)
 =\prod_{i=1}^{d}\prod_{j=1}^{b}
 \binom{h_{i,j-1}-h_{i+1,j}}{h_{i,j}-h_{i+1,j}}.
\end{equation}
\end{thm}
\pf
Given a tuple \(\boldsymbol\alpha\in\mathscr W_{k,d}\) satisfying
\(h_{\boldsymbol\alpha}=h\), condition \eqref{eq:sec4-14} assigns a subset \(S(C)\) to each
component \(C\in\mathcal B\). These subsets satisfy \eqref{eq:sec4-15}.

Conversely, let \(S:\mathcal B\to2^{\{0,\ldots,k-1\}}\) satisfy \eqref{eq:sec4-15}.
Extend \(S\) to cells by putting \(S(x)=S(C)\) for \(x\in C\),
\(S(x)=\{0,\ldots,k-1\}\) for \(x\in U_k\), and \(S(x)=\varnothing\) for \(x\notin U_1\).

For \(0\leq s<k\), define \(\alpha^{(s)} = \{x:s\in S(x)\}\). We claim that each \(\alpha^{(s)}\) is
a Young diagram. Suppose \(y\in\alpha^{(s)}\) and \(x\preceq y\). Then \(s\in S(y)\). By the
containment condition, \(S(y)\subseteq S(x)\), and therefore \(s\in S(x)\). Thus
\(x\in\alpha^{(s)}\). Hence \(\alpha^{(s)}\) is an order ideal.

Furthermore, every cell \(x\) belongs to exactly \(|S(x)|=h(x)\) component diagrams.
Thus the resulting tuple \(\boldsymbol\alpha\) satisfies \(h_{\boldsymbol\alpha}=h\).

The two constructions are inverse. By Theorem~\ref{thm:unique-representative}, the
tuples with multiplicity function \(h\) form one equivalence class. This proves the first
equality of \eqref{eq:sec4-20}.

For the product, we count the possible subsets
\(S_{i,j}=S(i,j)\subseteq\{0,\ldots,k-1\}\) directly.

Fill the cells row by row from bottom to top, and within each row from left to right.
Put \(S_{i,0}=\{0,\ldots,k-1\}\) and \(S_{d+1,j}=\varnothing\).
When \(S_{i,j}\) is to be chosen, both \(S_{i,j-1}\) and \(S_{i+1,j}\) have already been
chosen. The inclusion \(S_{i+1,j}\subseteq S_{i,j-1}\) already holds: for \(j>1\),
\(S_{i+1,j}\subseteq S_{i+1,j-1}\subseteq S_{i,j-1}\); for \(j=1\), it follows from
\(S_{i,0}=\{0,\ldots,k-1\}\). The Young-diagram conditions require
\(S_{i+1,j}\subseteq S_{i,j}\subseteq S_{i,j-1}\), with \(|S_{i,j}|=h_{i,j}\).
The subset \(S_{i+1,j}\) already provides \(h_{i+1,j}\) required elements. We therefore
choose \(h_{i,j}-h_{i+1,j}\) additional elements from \(S_{i,j-1}\setminus S_{i+1,j}\),
which has cardinality \(h_{i,j-1}-h_{i+1,j}\).

Hence the number of possible choices for \(S_{i,j}\) is
\[
\binom{
h_{i,j-1}-h_{i+1,j}
}{
h_{i,j}-h_{i+1,j}
}.
\]

Multiplying over all cells yields \eqref{eq:sec4-20}. \qed

When \(k=2\), the set \(\mathcal B\) consists of the edge-connected components of
\(U_1\setminus U_2\), and distinct components are incomparable under \(\succeq\). Hence
\begin{equation}
\label{eq:sec4-18}
 \Omega_2(U_1,U_2)=2^{\kappa(U_1\setminus U_2)},
\end{equation}
where \(\kappa\) is the number of edge-connected components. For \(k\geq3\), the
inclusions \(S(D)\subseteq S(C)\) required by \eqref{eq:sec4-15} generally prevent replacing
\(2^\kappa\) by \(k^\kappa\), even though counting the subsets one cell at a time gives the
product in \eqref{eq:sec4-20}.

\subsection{Generating functions for tuples and equivalence classes}
\label{subsec:plane-partitions}

Counting each cell by its multiplicity gives
\begin{equation}
\label{eq:sec4-22}
 \sum_{s=0}^{k-1}|\alpha^{(s)}|=\sum_x h(x)=\sum_{r=1}^{k}|U_r|.
\end{equation}
Thus every tuple in the equivalence class has size \(d^2+\sum_r|U_r|\).

\begin{thm}
\label{thm:tuple-class-generating-function}
We have
\begin{equation}
\label{eq:sec4-23}
 \sum_{d\geq0}\sum_{\substack{U_1\supseteq\cdots\supseteq U_k\\\ell(U_1)\leq d}}
 \Omega_k(U_1,\ldots,U_k)q^{d^2+\sum_r|U_r|}
 =\sum_{d\geq0}\frac{q^{d^2}}{(q;q)_d^k}.
\end{equation}
\end{thm}
\pf
By Theorem~\ref{thm:unique-representative}, every tuple belongs to a unique
equivalence class represented by \((d;U_1,\ldots,U_k)\), where
\(U_1\supseteq\cdots\supseteq U_k\). By Theorem~\ref{thm:class-cardinality}, the class contains
\(\Omega_k(U_1,\ldots,U_k)\) tuples, each of size \(d^2+\sum_{r=1}^{k}|U_r|\) by
\eqref{eq:sec4-22}. This gives the left-hand side of \eqref{eq:sec4-23}.

For the right-hand side, fix \(d\). Each of the \(k\) component partitions has at most \(d\)
parts and contributes \(1/(q;q)_d\). These partitions are chosen independently, and the
term \(d^2\) in the size contributes the factor \(q^{d^2}\). Thus the generating function is
\(q^{d^2}/(q;q)_d^k\). Summing over \(d\geq0\) proves the identity. \qed

For \(k=2\),
\[
\sum_{d\geq0}
\frac{q^{d^2}}{(q;q)_d^2}
=
\frac1{(q;q)_\infty},
\]
by the classical Durfee-square decomposition~\cite{and76}. Together with \eqref{eq:sec4-18}, Theorem~\ref{thm:tuple-class-generating-function} reduces exactly to the Boolean decomposition of Andrews and Dastidar~\cite{and26},
\[
\sum_{\lambda\in P^+}
2^{\kappa(\lambda)}q^{|\lambda|}
=
\frac1{(q;q)_\infty}.
\]

Here \(P^+\) denotes the partitions whose successive ranks are all nonnegative, and \(\kappa(\lambda)\) is the number of diagonal blocks.

We now count each equivalence class once, using its representative
\((d;U_1,\ldots,U_k)\). Thus the factor \(\Omega_k\) is omitted.

\begin{thm}
\label{thm:representatives-plane-partitions}

For fixed \(d\),
\begin{equation}
\label{eq:sec4-24}
\sum_{\substack{
U_1\supseteq\cdots\supseteq U_k\\
\ell(U_1)\leq d
}}
q^{\sum_{r=1}^{k}|U_r|}
=
\prod_{i=1}^{d}
\prod_{r=1}^{k}
\frac1{1-q^{i+r-1}}.
\end{equation}

Consequently,
\begin{equation}
\label{eq:sec4-25}
\sum_{d\geq0}
\sum_{\substack{
U_1\supseteq\cdots\supseteq U_k\\
\ell(U_1)\leq d
}}
q^{d^2+\sum_{r=1}^{k}|U_r|}
=
\sum_{d\geq0}
q^{d^2}
\prod_{i=1}^{d}
\prod_{r=1}^{k}
\frac1{1-q^{i+r-1}}.
\end{equation}

\end{thm}

\pf
Pad each \(U_r\) with zeros to \(d\) parts and form the array
\[
\Pi
=
\begin{pmatrix}
(U_1)_1&(U_1)_2&\cdots&(U_1)_d\\
(U_2)_1&(U_2)_2&\cdots&(U_2)_d\\
\vdots&\vdots&\ddots&\vdots\\
(U_k)_1&(U_k)_2&\cdots&(U_k)_d
\end{pmatrix}.
\]

Since each \(U_r\) is a partition, \((U_r)_1\geq (U_r)_2\geq\cdots\geq (U_r)_d\). Thus the rows of
\(\Pi\) are weakly decreasing.

The inclusion \(U_r\supseteq U_{r+1}\) means \((U_r)_i\geq(U_{r+1})_i\) for all \(i\), so the columns
are also weakly decreasing.

Hence \(\Pi\) is precisely a plane partition whose shape is contained in a \(k\times d\) rectangle.
Moreover, \(|\Pi| = \sum_{r=1}^{k}|U_r|\). Conversely, every plane partition whose shape
is contained in a \(k\times d\) rectangle determines a unique sequence
\(U_1\supseteq\cdots\supseteq U_k\) by reading its rows.

The classical MacMahon formula for plane partitions whose shape is contained in a
\(k\times d\) rectangle~\cite{mac16} gives
\[
\sum_{\Pi}q^{|\Pi|}
=
\prod_{r=1}^{k}
\prod_{i=1}^{d}
\frac1{1-q^{r+i-1}},
\]
which proves \eqref{eq:sec4-24}. Multiplying by \(q^{d^2}\) and summing over \(d\) gives \eqref{eq:sec4-25}. \qed

The two generating functions count different sets. In \eqref{eq:sec4-23}, the class
represented by \((d;U_1,\ldots,U_k)\) contributes
\(\Omega_k(U_1,\ldots,U_k)q^{d^2+\sum_r|U_r|}\), counting all of its tuples.
In \eqref{eq:sec4-25}, the same class contributes only \(q^{d^2+\sum_r|U_r|}\), counting its
representative once. For fixed \(d\), these representatives correspond to plane partitions
whose shape is contained in a \(k\times d\) rectangle.

\subsection{Cyclic shifts and pairing width}
\label{subsec:cyclic-pairing-connection}

Define the cyclic shift
\[
\sigma(d;\alpha^{(0)},\ldots,\alpha^{(k-1)})
 =(d;\alpha^{(1)},\ldots,\alpha^{(k-1)},\alpha^{(0)}).
\]
A tuple is fixed by \(\sigma\) if and only if all its component partitions are equal.
Equivalently, \(U_r\setminus U_{r+1}=\varnothing\) for \(1\leq r<k\), or
\(h_{\boldsymbol\alpha}(x)\in\{0,k\}\) for every cell \(x\).
Every fixed tuple has the form \((d;I,\ldots,I)\), where \(\ell(I)\leq d\), and has size
\(d^2+k|I|\). Their generating function is
\[
Z_k(q):=\sum_{\substack{\boldsymbol\alpha\in\mathscr W_k\\
 \sigma(\boldsymbol\alpha)=\boldsymbol\alpha}}q^{N(\boldsymbol\alpha)}
 =\sum_{d\geq0}\frac{q^{d^2}}{(q^k;q^k)_d}.
\]
For \(k=2\), tuples with equal components correspond to self-conjugate partitions, giving
\(Z_2(q)=\sum_{n\geq0}\operatorname{sc}(n)q^n\). For \(k>2\), \(Z_k\) counts tuples in
\(\mathscr W_k\) fixed by \(\sigma\), not ordinary self-conjugate partitions.

Put \(c_{\mathrm{col}}(\boldsymbol\alpha)=\sum_{s=0}^{k-1}s\alpha^{(s)}_1\).
Counting columns of a partition with at most \(d\) parts gives
\begin{equation}
\label{eq:sec4-31}
 \sum_{\ell(\alpha)\leq d}z^{\alpha_1}q^{|\alpha|}=\frac1{(zq;q)_d}.
\end{equation}

\begin{thm}
\label{thm:cyclic-root-of-unity}
Let \(\zeta\) be a primitive \(k\)-th root of unity, \(0\leq t<k\),
\(g=\gcd(k,t)\), and \(L=k/g\), with \(\gcd(k,0)=k\). Then
\begin{equation}
\label{eq:sec4-33}
 \sum_{\boldsymbol\alpha\in\mathscr W_k}
 \zeta^{t c_{\mathrm{col}}(\boldsymbol\alpha)}q^{N(\boldsymbol\alpha)}
 =\sum_{d\geq0}\frac{q^{d^2}}{(q^L;q^L)_d^g}
 =\sum_{\substack{\boldsymbol\alpha\in\mathscr W_k\\
 \sigma^t(\boldsymbol\alpha)=\boldsymbol\alpha}}q^{N(\boldsymbol\alpha)}.
\end{equation}
\end{thm}
\pf
Fix \(d\). By \eqref{eq:sec4-31}, the component with index \(s\) contributes
\(1/(\zeta^{ts}q;q)_d\). Thus the contribution for this value of \(d\) to the left-hand side
of \eqref{eq:sec4-33} equals
\[
q^{d^2}
\prod_{s=0}^{k-1}
\frac1{(\zeta^{ts}q;q)_d}.
\]

Expanding the finite \(q\)-Pochhammer symbols,
\begin{equation}
\label{eq:sec4-34}
q^{d^2}
\prod_{j=1}^{d}
\prod_{s=0}^{k-1}
\frac1{1-\zeta^{ts}q^j}.
\end{equation}

The root \(\zeta^t\) has order \(L\). As \(s\) runs through \(0,\ldots,k-1\), every \(L\)-th root of unity appears exactly \(g\) times. Hence
\[
\prod_{s=0}^{k-1}
(1-\zeta^{ts}x)
=
(1-x^L)^g.
\]

Substituting \(x=q^j\) into \eqref{eq:sec4-34} gives \(\frac{q^{d^2}} {(q^L;q^L)_d^g}\). This proves
the first equality.

Now consider the fixed points of \(\sigma^t\). Its action on the \(k\) component positions
has \(g\) cycles, each of length \(L\). A tuple is fixed precisely when its component
partitions are equal within each cycle. For fixed \(d\), it is therefore determined by \(g\)
arbitrary partitions \(I^{(1)},\ldots,I^{(g)}\), with \(\ell(I^{(j)})\leq d\). Each such partition is repeated \(L\) times, so the total size is \(d^2 +
L\sum_{j=1}^{g}|I^{(j)}|\). Its generating function is again \(\frac{q^{d^2}} {(q^L;q^L)_d^g}\).
Summing over \(d\) completes the proof. \qed

For fixed size \(n\), define
\(\mathscr W_k(n)=\{\boldsymbol\alpha:N(\boldsymbol\alpha)=n\}\) and
\(P_{k,n}(z)=\sum_{\boldsymbol\alpha\in\mathscr W_k(n)}z^{c_{\mathrm{col}}(\boldsymbol\alpha)}\).
Coefficient extraction in \eqref{eq:sec4-33} gives
\(P_{k,n}(\zeta^t)=\#\operatorname{Fix}_{\mathscr W_k(n)}(\sigma^t)\), the cyclic sieving
statement. Taking \(t=1\) in \eqref{eq:sec4-33} gives \(Z_k(q)\); equivalently,
\(P_{k,n}(\zeta)=[q^n]Z_k(q)\).

\begin{core}
\label{cor:pairing-width-conjugation}
For \(m\geq2\) and \(d\geq0\), the map \(\lambda\mapsto\alpha=(\Phi_m(\lambda))'\) is a
weight-preserving bijection
\[
\{\lambda:W_m(\lambda)\leq d\}\longleftrightarrow\{\alpha:\ell(\alpha)\leq d\},
\]
with \(\ell(\alpha)=W_m(\lambda)\) and \(\alpha_1=T_m(\lambda)\).
\end{core}
\pf
This follows from Corollary~\ref{cor:k-pairing-index-width}, since conjugation exchanges
length and largest part. \qed

Applying this bijection to each component identifies \(\mathscr W_{k,d}\) with the tuples
\((d;\lambda^{(0)},\ldots,\lambda^{(k-1)})\) satisfying
\(W_m(\lambda^{(s)})\leq d\) for every \(0\leq s<k\). Under this bijection,
\[
c_{\mathrm{col}}(\boldsymbol\alpha)=\sum_{s=0}^{k-1}sT_m(\lambda^{(s)}).
\]
The componentwise bijection commutes with cyclic shifts. Tuples with equal component
partitions correspond exactly to tuples \((d;\lambda,\ldots,\lambda)\). Consequently, the
statements on fixed points and root-of-unity evaluations also hold for the tuples satisfying
\(W_m(\lambda^{(s)})\leq d\) for every \(s\), with size
\(d^2+\sum_s|\lambda^{(s)}|\). In particular, the generating function for tuples with equal
components is
\[
\sum_{d\geq0}q^{d^2}\sum_{W_m(\lambda)\leq d}q^{k|\lambda|}=Z_k(q).
\]
For each component, pairing width and index are the numbers of rows and columns of the
corresponding Young diagram, respectively.

\vspace{0.5cm}
\noindent{\bf Acknowledgments.} This work was supported by  the National Science Foundation of China.

\end{document}